\documentclass[final]{elsarticle}
\usepackage{hyperref}
\usepackage[textwidth=16cm, textheight=24cm, left=2.5cm, right=2.5cm, top=3cm, bottom=3cm]{geometry}
\usepackage{latexsym}
\usepackage{enumerate}
\usepackage{amsmath, amsthm}
\usepackage{amssymb}
\usepackage{graphics}
\usepackage{graphicx}
\usepackage{subfigure}
\usepackage{indentfirst}
\usepackage{mathrsfs}
\usepackage{tikz}
\usepackage{longtable}
\usepackage{multirow}

\usetikzlibrary{decorations.pathreplacing}

\theoremstyle{plain}

\theoremstyle{definition}

\numberwithin{equation}{section}

\makeatletter
\def\ps@pprintTitle{%
  \let\@oddhead\@empty
  \let\@evenhead\@empty
  \def\@oddfoot{\reset@font\hfil\thepage\hfil}
  \let\@evenfoot\@oddfoot
}
\makeatother

\begin{document}
\begin{frontmatter}
\title{A general product of tensors with applications \tnoteref{title1}}
\tnotetext[title1]{Research supported by the National Science Foundation of China No.11231004}
\author[a]{Jia-Yu Shao\corref{cor1}}
\ead{jyshao@tongji.edu.cn}
\cortext[cor1]{Corresponding author.}

\address[a]{Department of Mathematics,  Tongji University,  Shanghai  200092,  China}

\date{2012-12-1}
\begin{abstract}

We define a general product of two $n$-dimensional tensors $\mathbb {A}$ and $\mathbb {B}$ with orders $m\ge 2$ and $k\ge 1$, respectively. This product is a generalization of the usual matrix product, and satisfies the associative law. Using this product, many concepts and known results of tensors can be simply expressed and/or proved, and a number of applications of this product will be given. Using this tensor product and some properties on the resultant of a system of homogeneous equations on $n$ variables, we define the similarity and congruence of tensors (which are also the generalizations of the corresponding relations for matrices), and prove that similar tensors have the same characteristic polynomials. We study two special kinds of similarity: permutational similarity and diagonal similarity, and their applications in the study of the spectra of hypergraphs and nonnegative irreducible tensors. We define the direct product of tensors (in matrix case it is also called the Kronecker product), and give its applications in the study of the spectra of two kinds of the products of hypergraphs. We also give applications of this tensor product in the study of nonnegative tensors, including a characterization of primitive tensors, the upper bounds of primitive degrees and the cyclic indices of some nonnegative irreducible tensors.

\end{abstract}

\end{frontmatter}

\section{Introduction}

In recent years, the study of tensors and the spectra of tensors (and hypergraphs) with their various applications has attracted much attention and interest, since the work of L.Qi ([12]) and L.H.Lim ([8]) in 2005.

\vskip 0.1cm

As is in [12], an order $m$ dimension $n$ tensor $\mathbb {A}=(a_{i_1i_2\cdots i_m})_{1\le i_j\le n \ (j=1,\cdots ,m)}$ over the complex field $\mathbb {C}$ is a multidimensional array with all entries $a_{i_1i_2\cdots i_m}\in \mathbb {C} \  \ (i_1,\cdots ,i_m\in [n]=\{1,\cdots ,n\})$. For a vector $x=(x_1,\cdots ,x_n)^T\in \mathbb {C}^n$, let $\mathbb {A}x^{m-1}$ be a vector in $\mathbb {C}^n$ whose $i$th component is defined as the following:
$$(\mathbb {A}x^{m-1})_i=\sum_{i_2,\cdots ,i_m=1}^n a_{ii_2\cdots i_m}x_{i_2}\cdots x_{i_m} \eqno {(1.1)}$$
and let $x^{[r]}=(x_1^r,\cdots ,x_n^r)^T$. Then ([3,12]) a number $\lambda \in \mathbb {C}$ is called an eigenvalue of $\mathbb {A}$ if there exists a nonzero vector $x\in \mathbb {C}^n$ such that
$$\mathbb {A}x^{m-1}=\lambda x^{[m-1]} \eqno {(1.2)}$$
and in this case, $x$ is called an eigenvector of $\mathbb {A}$ corresponding to the eigenvalue $\lambda$. L.Qi also defined several other types of eigenvalues (and eigenvectors) in [12].

\vskip 0.1cm

Since L.Qi introduced the definitions of various types of the eigenvalues (and eigenvectors) of tensors, quite many interesting research works in this area appeared (see [2-5] and [8-15]). But noticing that, the term $\mathbb {A}x^{m-1}$ in the left side of (1.1) is only a special notation (its inside part $x^{m-1}$ does not have specific meaning), and also we can not view the notation $\mathbb {A}x^{m-1}$ as some general operation between $\mathbb {A}$ and $x^{m-1}$. Anyway, we have not seen some kinds of general product between tensors up to now.

\vskip 0.1cm

In this paper, we will introduce (in Definition 1.1) a general tensor product defined for any two tensors of dimension $n$ (one has order $m\ge 2$ and another one has order $k\ge 1$). We will see that, the product of two dimension $n$ tensors is still a tensor of dimension $n$ (but in general will have different order). This tensor product also possess many properties. Among them the most important and useful property is the associative law (see Theorem 1.1).

\vskip 0.1cm

Besides these properties, we can also see that, when the two tensors are matrices (or the first one is a matrix and the second one is a vector), our definition of tensor product coincides with the usual matrix product, so it is a generalization of the matrix product.

\vskip 0.1cm

In this paper, we will also give a number of applications of this tensor product. we will see that, by using this tensor product, not only many mathematical expressions can be simplified (for example, the term $\mathbb {A}x^{m-1}$ in the left side of (1.1) can now be directly simplified as $\mathbb {A}x$), and many known results can be simply stated and proved, but we can also obtain many new results (see sections 2,3 and 4).

\vskip 0.1cm

In section 2, we first use this tensor product and some properties on resultant (namely, the hyperdeterminant)
to define the similarity of tensors. We prove that the similar tensors have the same characteristic polynomials, thus have the same spectrum (as multiset). We also study two special kinds of similarity: permutation similarity and diagonal similarity. Using the permutation similarity we can show that the spectrum of a hypergraph is independent of the ordering of its vertices (but the adjacency tensor of a hypergraph depends on the ordering of its vertices). Using the diagonal similarity, we obtain a result on the symmetry of the spectrum of a nonnegative irreducible tensor which is stronger than the result in [14]. That is: if a nonnegative irreducible tensor $\mathbb {A}$ has cyclic index $k$, then the whole spectrum of $\mathbb {A}$ is invariant under a rotation of angle $2\pi /k$ in the complex plane. We also define in section 2 the congruence relation between tensors and apply this relation to the study of {\bf E}-eigenvalues and {\bf E}-eigenvectors ([12]).

\vskip 0.1cm

In section 3, we generalize the direct product of matrices to the direct product of tensors (of the same order, but may be different dimensions). We first prove a relation (Theorem 3.1) between the direct product of tensors and the general product of tensors defined in Definition 1.1. Then we use this to study the spectra of two kinds of products of hypergraphs.

\vskip 0.1cm

In section 4, we apply this tensor product to the study of nonnegative tensors. We give a simple characterization for the primitivity of a nonnegative tensor in terms of the powers of $\mathbb {A}$. We show that the primitivity of a nonnegative tensor $\mathbb {A}$ is uniquely determined by the zero pattern of $\mathbb {A}$. From this we show the boundedness of the primitive degrees (for fixed $n$ and $m$), and give some upper bounds for the primitive degrees of the primitive tensors. We also propose a conjecture concerning the upper bound of the primitive degrees.

\vskip 0.28cm

For the sake of simplicity, we sometime use the following "binary notation" for the subscripts of the tensor. Namely, we will write $a_{i_1i_2\cdots i_m}$ as $a_{i_1\alpha }$, where $\alpha =i_2\cdots i_m \in [n]^{m-1}$.

\vskip 0.1cm

Now we define the general product of two $n$-dimensional tensors.

\vskip 0.1cm

\noindent {\bf Definition 1.1}: Let $\mathbb {A}$ (and $\mathbb {B}$) be an order $m\ge 2$ (and order $k\ge 1$), dimension $n$ tensor, respectively. Define the product $\mathbb {A}\mathbb {B}$ to be the following tensor $\mathbb {C}$ of order $(m-1)(k-1)+1$ and dimension $n$:
$$c_{i\alpha_1\cdots \alpha_{m-1}}=\sum_{i_2,\cdots ,i_m=1}^n a_{ii_2\cdots i_m}b_{i_2\alpha_1}\cdots b_{i_m\alpha_{m-1}} \quad  (i\in [n], \  \alpha_1, \cdots ,\alpha_{m-1}\in [n]^{k-1})$$

\vskip 0.1cm

\noindent {\bf Example 1.1}: When $k=1$ and $\mathbb {B}=x\in \mathbb {C}^n$ is a vector of dimension $n$, then
$(m-1)(k-1)+1=1$. Thus $\mathbb {A}\mathbb {B}=\mathbb {A}x$ is still a vector of dimension $n$, and we have
$$(\mathbb {A}x)_i=(\mathbb {A}\mathbb {B})_i=c_i=\sum_{i_2,\cdots ,i_m=1}^n a_{ii_2\cdots i_m}x_{i_2}\cdots x_{i_m}=(\mathbb {A}x^{m-1})_i  \quad (\mbox {where } \ \mathbb {A}x^{m-1} \  \mbox {is defined by (1.2)})$$
Thus we have $\mathbb {A}x^{m-1}=\mathbb {A}x$. So the first application of the tensor product defined above is that now $\mathbb {A}x^{m-1}$ can be simply written as $\mathbb {A}x$.

\vskip 0.1cm

\noindent {\bf Example 1.2}: (1). If $k=2$ (namely $\mathbb {B}$ is a matrix of order $n$), then $(m-1)(k-1)+1=m$. So in this case $\mathbb {A}\mathbb {B}$ is still a tensor of order $m$.

\vskip 0.1cm

\noindent (2). If $m=2$ (namely $\mathbb {A}$ is a matrix of order $n$), then $(m-1)(k-1)+1=k$. So in this case $\mathbb {A}\mathbb {B}$ is still a tensor of order $k$.

\vskip 0.1cm

\noindent (3). It is easy to check from the definition that $I_n\mathbb {A}=\mathbb {A}=\mathbb {A}I_n$, where $I_n$ is the identity matrix of order $n$.

\vskip 0.1cm

\noindent {\bf Example 1.3}: In case when both $\mathbb {A}$ and $\mathbb {B}$ are matrices, or when $\mathbb {A}$ is a matrix and $\mathbb {B}$ is a vector, our tensor product $\mathbb {A}\mathbb {B}$ coincides with the usual matrix product. So it is a generalization of the matrix product.

\vskip 0.38cm

\noindent {\bf Proposition 1.1}: The tensor product defined in Definition 1.1 has the following properties.

\vskip 0.1cm

\noindent (1). (The left distributivity): \quad  $(\mathbb {A}_1+\mathbb {A}_2)\mathbb {B}=\mathbb {A}_1\mathbb {B}+\mathbb {A}_2\mathbb {B}.$

\vskip 0.1cm

\noindent (2). (The right distributivity in the case of $m=2$): \quad  $\mathbb {A}(\mathbb {B}_1+\mathbb {B}_2)=\mathbb {A}\mathbb {B}_1+\mathbb {A}\mathbb {B}_2 \ $
(when $\mathbb {A}$ is a matrix).

\vskip 0.1cm

\noindent (Note that in general when $\mathbb {A}$ is not a matrix, then the right distributivity does not hold.)

\vskip 0.1cm

\noindent (3).  $(\lambda \mathbb {A})\mathbb {B}=\lambda (\mathbb {A}\mathbb {B}).$  \qquad ($\lambda \in \mathbb {C}$.)

\vskip 0.1cm

\noindent (4). $\mathbb {A}(\lambda \mathbb {B})=\lambda ^{m-1}(\mathbb {A}\mathbb {B}).$  \qquad ($\lambda \in \mathbb {C}$.)

\vskip 0.2cm

Next we will show that the associative law of this tensor product holds. The following identity will be needed in the proof of the associative law.

$$\prod_{j=1}^m\sum_{t_{j1},\cdots ,t_{jk}=1}^n f(j,t_{j1},\cdots ,t_{jk})= \sum_{t_{jh}=1 (1\le j\le m; 1\le h\le k)}^n\prod_{j=1}^m f(j,t_{j1},\cdots ,t_{jk})  \eqno {(1.1)}$$

For the simplicity of the notation, we assume in the following that the tensors $\mathbb {A}$ ,$\mathbb {B}$ and $\mathbb {C}$ have the orders $m+1$, $k+1$ and $r+1$, respectively.

\vskip 0.3cm

\noindent {\bf Theorem 1.1 (The associative law of the tensor product)}: Let $\mathbb {A}$ (and $\mathbb {B}$, $\mathbb {C}$) be an order $m+1$ (and order $k+1$, order $r+1$), dimension $n$ tensor, respectively.  Then we have
$$\mathbb {A}(\mathbb {B}\mathbb {C})=(\mathbb {A}\mathbb {B})\mathbb {C}$$
\noindent {\bf Proof}. For $\beta_1,\cdots ,\beta_m\in ([n]^r)^k$, we write:
$$\beta_1=\theta_{11}\cdots \theta_{1k}, \  \  \cdots ,  \   \  \beta_m=\theta_{m1}\cdots \theta_{mk} \qquad (\theta_{ij}\in [n]^r, \ i=1,\cdots ,m; \ j=1,\cdots ,k.) $$
Then we have:
$$
\begin{aligned}
\displaystyle & (\mathbb {A}(\mathbb {B}\mathbb {C}))_{i\beta_1\cdots \beta_m}= \sum_{i_1,\cdots ,i_m=1}^n a_{ii_1\cdots i_m}\left (\prod_{j=1}^m (\mathbb {B}\mathbb {C})_{i_j\beta_j} \right )        \\
&= \sum_{i_1,\cdots ,i_m=1}^n a_{ii_1\cdots i_m}\left (\prod_{j=1}^m (\mathbb {B}\mathbb {C})_{i_j\theta_{j1}\cdots \theta_{jk}}     \right )          \\
&= \sum_{i_1,\cdots ,i_m=1}^n a_{ii_1\cdots i_m}\left (\prod_{j=1}^m \sum_{t_{j1},\cdots ,t_{jk}=1}^n b_{i_jt_{j1}\cdots t_{jk}} (c_{t_{j1}\theta_{j1}}\cdots c_{t_{jk}\theta_{jk}} )    \right )         \\
&\stackrel{(*)}{=} \sum_{i_1,\cdots ,i_m=1}^n a_{ii_1\cdots i_m} \sum_{t_{jh}=1 (1\le j\le m; 1\le h\le k)}^n  \left ( \prod_{j=1}^m   b_{i_jt_{j1}\cdots t_{jk}} (c_{t_{j1}\theta_{j1}}\cdots c_{t_{jk}\theta_{jk}} )    \right )        \\
\end{aligned} \eqno  {(1.2)} $$
where the equation (*) follows from equation (1.1).

\vskip 0.1cm

On the other hand, for $\alpha_1,\cdots ,\alpha_m\in [n]^k$, we write:
$$\alpha_1=t_{11}\cdots t_{1k}, \  \   \cdots ,  \   \  \alpha_m=t_{m1}\cdots t_{mk} \qquad (t_{ij}\in [n], \ i=1,\cdots ,m; \ j=1,\cdots ,k.) $$
Then we also have:
$$
\begin{aligned}
\displaystyle & ((\mathbb {A}\mathbb {B})\mathbb {C})_{i\beta_1\cdots \beta_m}= \sum_{\alpha_1,\cdots ,\alpha_m\in [n]^k} (\mathbb {A}\mathbb {B})_{i\alpha_1\cdots \alpha_m}\left (\prod_{j=1}^m (c_{t_{j1}\theta_{j1}}\cdots c_{t_{jk}\theta_{jk}} ) \right )        \\
&= \sum_{t_{jh}=1 (1\le j\le m; 1\le h\le k)}^n \sum_{i_1,\cdots ,i_m=1}^n a_{ii_1\cdots i_m} \left (\prod_{j=1}^m b_{{i_j}\alpha_j}  \right )   \left (\prod_{j=1}^m (c_{t_{j1}\theta_{j1}}\cdots c_{t_{jk}\theta_{jk}} ) \right )     \\
&=  \sum_{i_1,\cdots ,i_m=1}^n a_{ii_1\cdots i_m} \sum_{t_{jh}=1 (1\le j\le m; 1\le h\le k)}^n  \left ( \prod_{j=1}^m   b_{i_jt_{j1}\cdots t_{jk}} (c_{t_{j1}\theta_{j1}}\cdots c_{t_{jk}\theta_{jk}} )    \right )          \\
\end{aligned} \eqno  {(1.3)} $$
Comparing the right hand sides of (1.2) and (1.3), we obtain $\mathbb {A}(\mathbb {B}\mathbb {C})=(\mathbb {A}\mathbb {B})\mathbb {C}$.
\qed

\vskip 0.2cm

By the associative law, we can henceforth write $(\mathbb {A}\mathbb {B})\mathbb {C}$ and $\mathbb {A}(\mathbb {B}\mathbb {C})$ as $\mathbb {A}\mathbb {B}\mathbb {C}$. We can also define $\mathbb {A}^k$ as the product of $k$ many tensors $\mathbb {A}$. Furthermore, by the left distributive law and right distributive law for matrices, we have:
$$P(\mathbb {A}+\mathbb {B})Q=P\mathbb {A}Q+P\mathbb {B}Q \  \qquad  \mbox {(when $P,Q$ are both matrices)} \eqno {(1.4)} $$

\vskip 0.2cm

The unit tensor of order $m$ and dimension $n$ is the tensor $\mathbb {I}=(\delta _{i_1,i_2,\cdots ,i_m})$ with entries as follows:
$$\delta _{i_1,i_2,\cdots ,i_m}=\left\{
                            \begin{array}{cc}
                              1 &   \mbox{if } i_1=i_2=\cdots =i_m  \\
                              0 &  \mbox{otherwise}\\
                            \end{array}
                          \right.$$
It is easy to see from the definition that $\mathbb {I}x=x^{[m-1]}$. Thus equation (1.2) can be rewritten as
$$(\lambda \mathbb {I}- \mathbb {A})x=0  \eqno {(1.5)}$$

\vskip 1.88cm

\section{The similarity and congruence of tensors}

In this section, we first use the tensor product defined in section 1, together with some properties on the resultant of a system of homogeneous equations on $n$ variables, to define the similarity of tensors (which is a generalization of the similarity of matrices), and then prove that the similar tensors have the same characteristic polynomials, hence the same spectrum (as a multiset). We further study two special kinds of similarity: permutation similarity and diagonal similarity, then apply them respectively in the study of the spectra of hypergraphs, and in the study of the symmetry of the spectra of nonnegative irreducible tensors. We also define the congruence relation between tensors, and show how this relation can be applied in the study of {\bf E}-eigenvalues and {\bf E}-eigenvectors of tensors ([12]).

In [6], D.Cox et al used the resultant to define the hyperdeterminant $det(\mathbb {A})$ of a tensor $\mathbb {A}$.

\vskip 0.18cm

\noindent {\bf Definition 2.1} ([6], hyperdeterminant): Let $\mathbb {A}$ be an order $m$ dimension $n$ tensor. Then the hyperdeterminant of $\mathbb {A}$, denoted by $det(\mathbb {A})$, is the resultant of the system of homogeneous equations $\mathbb {A}x=0$, which is also the unique polynomial on the entries of $\mathbb {A}$ (viewed as independent variables) satisfying the following three conditions:

\vskip 0.1cm

\noindent (1) $det(\mathbb {A})=0$ if and only if the system of homogeneous equations $\mathbb {A}x=0$ has a nonzero solution.

\vskip 0.1cm

\noindent (2) $det(\mathbb {I})=1$, where $\mathbb {I}$ is the unit tensor.

\vskip 0.1cm

\noindent (3) $det(\mathbb {A})$ is an irreducible polynomial on the entries of $\mathbb {A}$, when the entries of $\mathbb {A}$ are viewed as independent variables.

\vskip 0.1cm

The proof of the existence and uniqueness of the polynomial satisfying the above three conditions can be found in [7].

\vskip 0.18cm

\noindent {\bf Definition 2.2}: Let $\mathbb {A}$ be an order $m$ dimension $n$ tensor. Then the characteristic polynomial of $\mathbb {A}$, denoted by $\phi_{\mathbb {A}}(\lambda )$, is the hyperdeterminant $det(\lambda \mathbb {I}- \mathbb {A})$.

\vskip 0.1cm

It is easy to see from the above discussions and the equation (1.5) that, $\lambda _0$ is an eigenvalue of $\mathbb {A}$ if and only if it is a root of the characteristic polynomial of $\mathbb {A}$.

\vskip 0.28cm

In the following study of the similarity and congruence of tensors, we need to use the following formula for the entries of the tensor $P\mathbb {A}Q$ (where $P,Q$ are both matrices).

$$(P\mathbb {A}Q)_{i_1\cdots i_m}=\sum_{j_1,\cdots ,j_m=1}^n a_{j_1\cdots j_m}p_{i_1j_1}q_{j_2i_2}\cdots q_{j_mi_m}  \eqno {(2.1)} $$
\vskip 0.1cm

By using the three conditions of the hyperdeterminant in Definition 2.1, we can prove the following lemma.

\vskip 0.28cm

\noindent {\bf Lemma 2.1}: Let $\mathbb {A}$ be an order $m$ dimension $n$ tensor,  $\mathbb {I}$ be the order $m$ dimension $n$ unit tensor, $P$ and $Q$ are two matrices of order $n$. Then we have:
$$det(P\mathbb {A}Q)=det(P\mathbb {I}Q)det(\mathbb {A})  \eqno {(2.2)}$$

\noindent {\bf Proof}. We consider the following two cases.
\vskip 0.1cm

\noindent {\bf Case 1}: Both $P,Q$ are invertible.

Then by the associative law of the tensor product, we can see that $\mathbb {A}x=0$ has a nonzero solution $x$ if and only if $(P\mathbb {A}Q)y=0$ has a nonzero solution $y$, which implies that $det(\mathbb {A})=0$ if and only if $det(P\mathbb {A}Q)=0$.

Now by the definition, $det(\mathbb {A})$ is irreducible, so $det(\mathbb {A})$ is a factor of $det(P\mathbb {A}Q)$. Similarly, $det(P\mathbb {A}Q)$ is a factor of $det(\mathbb {A})$ since $P,Q$ are both invertible. Thus there exists some constant $c=c(P,Q)$ (independent of $\mathbb {A}$) such that $det(P\mathbb {A}Q)=c\cdot  det(\mathbb {A})$.

Since $c$ is independent of $\mathbb {A}$, we can take $\mathbb {A}=\mathbb {I}$, which gives us that $c=det(P\mathbb {I}Q)$, so (2.2) holds in this case.

\noindent {\bf Case 2}: The general case.

Take $P_{\epsilon }=P+\epsilon I$ and $Q_{\epsilon }=Q+\epsilon I$. Then $P_{\epsilon }$ and $Q_{\epsilon }$ can be both invertible when $\epsilon $ is in some small open interval $(0,\delta )$. Thus by Case 1 we have
$$det(P_{\epsilon }\mathbb {A}Q_{\epsilon })=det(P_{\epsilon }\mathbb {I}Q_{\epsilon })det(\mathbb {A})
 \qquad  (\epsilon \in (0,\delta ))   \eqno {(2.3)}$$
Take the limit $\epsilon \rightarrow 0$ on both sides of (2.3), we obtain (2.2).
\qed

\vskip 0.2cm

\noindent {\bf Definition 2.3}: Let $\mathbb {A}$ and $\mathbb {B}$ be two order $m$ dimension $n$ tensors.  Suppose that there exist two matrices $P$ and $Q$  of order $n$ with $P\mathbb {I}Q =\mathbb {I}$ such that $\mathbb {B}=P\mathbb {A}Q$, then we say that the two tensors $\mathbb {A}$ and $\mathbb {B}$ are similar.

\vskip 0.28cm

\noindent {\bf Remark 2.1}: If $P$ and $Q$ are two matrices of order $n$ with $P\mathbb {I}Q =\mathbb {I}$, where $\mathbb {I}$ is the order $m$ dimension $n$ unit tensor, then both $P$ and $Q$ are invertible matrices.

\noindent {\bf Proof. (1).} If $Q$ is not invertible, then there exists some nonzero vector $x\in \mathbb {C}^n$ such that $Qx=0$. Thus by the associative law of the tensor product we have
$x^{[m-1]}=\mathbb {I}x=P\mathbb {I}Q x=0$, which implies that $x=0$, a contradiction.

\noindent {\bf (2).} If $P$ is not invertible, then there exists some nonzero vector $y\in \mathbb {C}^n$ such that $Py=0$. Now take $x=Q^{-1}y^{[\frac {1} {m-1}]}$, then $x$ is a nonzero vector since $y$ is. But then we would have
$x^{[m-1]}=\mathbb {I}x=P\mathbb {I}Q x=P\mathbb {I}y^{[\frac {1} {m-1}]}=Py=0$, which implies that $x=0$, a contradiction.
\qed

\vskip 0.28cm

\noindent {\bf Remark 2.2}: The similarity relation of tensors is an equivalent relation.

\noindent {\bf Proof. } The reflection follows from $I_n\mathbb {A}I_n=\mathbb {A}$ (where $I_n$ is the identity matrix). The symmetry follows from Remark 2.1 and the associative law of the tensor product. The transitivity also follows from the associative law of the tensor product.

\vskip 0.2cm

The following theorem shows that similar tensors have the same characteristic polynomials, and thus they have the same spectrum (as a multiset).

\vskip 0.2cm

\noindent {\bf Theorem 2.1}: Let $\mathbb {A}$ and $\mathbb {B}$ be two order $m$ dimension $n$ similar tensors.   Then we have:

\noindent {\bf (1).}   $det (\mathbb {B})=det (\mathbb {A})$.

\noindent {\bf (2).}  $\phi _{\mathbb {B}}(\lambda )=\phi _{\mathbb {A}}(\lambda )$.

\noindent {\bf Proof.} By definition, there exist two matrices $P$ and $Q$ of order $n$ with $P\mathbb {I}Q =\mathbb {I}$ and $\mathbb {B}=P\mathbb {A}Q$. Then (1) follows directly from Lemma 2.1.

For (2), by using the equation (1.4) and the relations $P\mathbb {I}Q =\mathbb {I}$ and $\mathbb {B}=P\mathbb {A}Q$, we have
$$P(\lambda \mathbb {I}-\mathbb {A})Q=\lambda P\mathbb {I}Q - P\mathbb {A}Q= \lambda \mathbb {I}-\mathbb {B} $$
Thus $\lambda \mathbb {I}-\mathbb {A}$ and $\lambda \mathbb {I}-\mathbb {B}$ are also similar.  Using this and the result (1) of this theorem we have
$$\phi _{\mathbb {B}}(\lambda )=det (\lambda \mathbb {I}-\mathbb {B})
=det (\lambda \mathbb {I}-\mathbb {A})=\phi _{\mathbb {A}}(\lambda ).$$
\qed

In the following, we consider two special kinds of the similarity of tensors. The first one is the permutation similarity of tensors, and the second one is the diagonal similarity of tensors (also see [14,15]).

\vskip 0.2cm

\noindent {\bf Lemma 2.2}:  Let $\sigma \in S_n$ be a permutation on the set [$n$], $P=P_{\sigma }=(p_{ij})$ be the corresponding permutation matrix of $\sigma $ (where $p_{ij}=1\Longleftrightarrow j=\sigma (i)$). Let $\mathbb {A}$ and $\mathbb {B}$ be two order $m$ dimension $n$ tensors with $\mathbb {B}=P\mathbb {A}P^T$. Then we have:

\noindent {\bf (1).}   $b_{i_1\cdots i_m}=a_{\sigma (i_1)\cdots \sigma (i_m)}$.

\noindent {\bf (2).}  $P\mathbb {I}P^T =\mathbb {I}$. \qquad (Thus $\mathbb {A}$ and $\mathbb {B}$ are similar.)

\noindent {\bf Proof}. (1) By using the formula (2.1)  we have:

$$b_{i_1\cdots i_m}=\sum _{j_1,\cdots ,j_m=1}^n a_{j_1\cdots j_m}p_{i_1j_1}((P^T)_{j_2i_2}\cdots (P^T)_{j_mi_m})=\sum _{j_1,\cdots ,j_m=1}^n a_{j_1\cdots j_m}p_{i_1j_1}p_{i_2j_2}\cdots p_{i_mj_m}=a_{\sigma (i_1)\cdots \sigma (i_m)}$$

\noindent {\bf (2).}  From (1) we see that the diagonal elements of $\mathbb {B}$ are those diagonal elements of $\mathbb {A}$, and the non-diagonal elements of $\mathbb {B}$ are those non-diagonal elements of $\mathbb {A}$. Thus we have $\mathbb {B}=\mathbb {I}$ when $\mathbb {A}=\mathbb {I}$. Namely $P\mathbb {I}P^T =\mathbb {I}$.
\qed

\vskip 0.2cm

From (2) of Lemma 2.2, we see that if $P=P_{\sigma }$ is a permutation matrix, $\mathbb {A}$ and $\mathbb {B}$ are two order $m$ dimension $n$ tensors with $\mathbb {B}=P\mathbb {A}P^T$, then $\mathbb {A}$ and $\mathbb {B}$ are similar tensors. In this case we say that $\mathbb {A}$ and $\mathbb {B}$ are {\bf permutational similar}.

\vskip 0.2cm

It is easy to see that if $\mathbb {A}$ and $\mathbb {B}$ are {\bf permutational similar}, and $\mathbb {A}$ is supersymmetric, then $\mathbb {B}$ is also supersymmetric.

\vskip 0.2cm

An important application of the permutational similar tensors is on the spectra of hypergraphs. According to the definition, the spectrum of a hypergraph is the spectrum of its adjacency tensor. But the adjacency tensor of a hypergraph depends on the ordering of its vertices. Thus a natural question arises: Is the spectrum of a hypergraph independent of the ordering of its vertices?

\vskip 0.2cm

By using the permutational similar tensors, it is easy to see that if we have two different orderings of the vertices of a hypergraph, then the two adjacency tensors are permutational similar. Since permutational similar tensors have the same spectra (by Theorem 2.1), we can now obtain an affirmative answer to the above question:

\vskip 0.2cm

\noindent {\bf Theorem 2.2}: (1).  The spectrum (as a multiset) of a hypergraph is independent of the ordering of its vertices.

\vskip 0.18cm

\noindent (2). Isomorphic hypergraphs have the same spectra.
\qed

\vskip 0.28cm

Now we consider the second special kind of similarity of tensors: the diagonal similarity. This relation of tensors was first studied in [14] and [15], and was applied to the study of the spectra of nonnegative irreducible tensors. Here we use our tensor product to show that, this relation between tensors is also a special kind of similarity in the sense of Definition 2.3.

\vskip 0.3cm

\noindent {\bf Lemma 2.3}: Let $D=diag (d_{11},\cdots , d_{nn})$ be an invertible diagonal matrix of order $n$. Then we have $D^{-(m-1)}\mathbb {I}D = \mathbb {I}$.

\noindent {\bf Proof. } By using the formula (2.1) we have
$$(D^{-(m-1)}\mathbb {I}D)_{i_1\cdots i_m}= \sum _{j_1,\cdots ,j_m=1}^n \delta _{j_1\cdots j_m}(D^{-(m-1)})_{i_1j_1}d_{j_2i_2}\cdots d_{j_mi_m}=\sum _{j=1}^n (D^{-(m-1)})_{i_1j}d_{ji_2}\cdots d_{ji_m}=\delta _{i_1\cdots i_m} $$
\qed

\vskip 0.2cm

\noindent {\bf Definition 2.4}: Let $\mathbb {A}$ and $\mathbb {B}$ be two order $m$ dimension $n$ tensors. We say that $\mathbb {A}$ and $\mathbb {B}$ are diagonal similar, if there exists some invertible diagonal matrix $D$ of order $n$ such that $\mathbb {B}=D^{-(m-1)}\mathbb {A}D$.

\vskip 0.3cm

\noindent {\bf Theorem 2.3}: If the two order $m$ dimension $n$ tensors $\mathbb {A}$ and $\mathbb {B}$ are diagonal similar, Then $\phi _{\mathbb {B}}(\lambda )=\phi _{\mathbb {A}}(\lambda )$.

\noindent {\bf Proof. }  From Lemma 2.3 we see that diagonal similar tensors are similar.  \qed

\vskip 0.2cm

Using Theorem 2.3, we can obtain the following result (which is stronger than a related result in [14]) on the symmetry of the spectrum of a nonnegative irreducible tensor.

\vskip 0.2cm

\noindent {\bf Theorem 2.4}: Let $\mathbb {A}$ be an order $m$ dimension $n$ nonnegative irreducible tensor with the spectral radius $\rho (A)$. Suppose that the cyclic index of $\mathbb {A}$ is $k$ (namely, $\mathbb {A}$ has exactly $k$ distinct eigenvalues with modulus $\rho (A)$). Then $spec (\mathbb {A})$ (as a multiset) is invariant under a rotation of an angle $\frac {2\pi } {k}$ on the complex plane, and therefore $\phi _{\mathbb {A}}(\lambda )=\lambda ^rf(\lambda ^k)$ for some nonnegative integer $r$ and some polynomial $f(x)$.

\vskip 0.2cm

\noindent {\bf Proof. } By [14], there exists an invertible diagonal matrix $D$, such that
$$\mathbb {A}=e^{-\frac {2\pi i} {k}}D^{-(m-1)}\mathbb {A}D.$$
So from Theorem 2.2 we have $spec (\mathbb {A})=e^{-\frac {2\pi i} {k}}spec (D^{-(m-1)}\mathbb {A}D)=e^{-\frac {2\pi i} {k}}spec (\mathbb {A}).$

Now suppose that
$$\phi _{\mathbb {A}}(\lambda )= \sum_{j=0}^da_j\lambda ^{d-j}, \quad  (\mbox {where $d$ is the degree of} \  \phi _{\mathbb {A}}(\lambda ))$$
Then from the above result we have
$$e^{\frac {2\pi di} {k}}\phi _{\mathbb {A}}(\lambda )=\phi _{\mathbb {A}}(e^{\frac {2\pi i} {k}}\lambda ), \qquad  \mbox {and thus}  \qquad
e^{\frac {2\pi di} {k}}\sum_{j=0}^da_j\lambda ^{d-j} = \sum_{j=0}^da_j(e^{\frac {2\pi i} {k}}\lambda )^{d-j}$$
Comparing the coefficients of both sides, we obtain $e^{\frac {2\pi di} {k}}a_j=e^{\frac {2\pi (d-j)i} {k}}a_j$. Thus $a_j=0$ if $j$ is not a multiple of $k$. Therefore we have $\phi _{\mathbb {A}}(\lambda )=\lambda ^rf(\lambda ^k)$.
\qed

\vskip 0.28cm

In [14] and [15], Q.Yang and Y.Yang give the following results regarding the diagonal similar tensors. Here we show that the proofs of these results can be simplified by using the associativity of our tensor product.

\vskip 0.2cm

\noindent {\bf Theorem 2.5 ([14] and [15]) }: (1). Suppose that the two tensors $\mathbb {A}$ and $\mathbb {B}$ are diagonal similar, namely $\mathbb {B}=D^{-(m-1)}\mathbb {A}D$ for some invertible diagonal matrix $D$. Then $x$ is an eigenvector of $\mathbb {B}$ corresponding to the eigenvalue $\lambda $ if and only if $y=Dx$ is an eigenvector of $\mathbb {A}$ corresponding to the same eigenvalue $\lambda $.

\vskip 0.2cm

\noindent (2). Every nonnegative irreducible tensor $\mathbb {A}$ is diagonal similar to a multiple of some stochastic tensor.

\vskip 0.2cm

\noindent {\bf Proof }: (1). By the associative law of the tensor product, we have for $y=Dx$ that
$$\mathbb {B}x=\lambda x^{[m-1]}\Longleftrightarrow D^{-(m-1)}\mathbb {A}Dx=\lambda x^{[m-1]}\Longleftrightarrow \mathbb {A}y=D^{(m-1)}\lambda x^{[m-1]}=\lambda (Dx)^{[m-1]}\Longleftrightarrow \mathbb {A}y=\lambda y^{[m-1]}$$

\noindent (2). Let $u$ be the positive eigenvector of $\mathbb {A}$ corresponding to the eigenvalue $\rho (\mathbb {A})$. Let $e$ be the column vector of dimension $n$ with all components 1, let $D$ be the invertible diagonal matrix with $u=De$, and let $\mathbb {B}=D^{-(m-1)}\mathbb {A}D$. Then $\mathbb {B}$ is also nonnegative. So from the result (1) and $\mathbb {A}u= \rho (\mathbb {A}) u^{[m-1]}$, we have $\mathbb {B}e=\rho (\mathbb {A}) e^{[m-1]}=\rho (\mathbb {A}) e$, which means that $\frac {1} {\rho (\mathbb {A})} \mathbb {B}$ is a stochastic tensor.           \qed

\vskip 0.58cm

Now we consider another relation between tensors: the congruence relation (which is also a generalization of the congruence relation between matrices) and its applications on the {\bf E}-eigenvalues and {\bf E}-eigenvectors
 of tensors. First we recall that in [12, section 6], L.Qi defined a useful relation between tensors (which he called "orthogonal similarity") as follows.

\vskip 0.2cm

\noindent {\bf Definition 2.5} (L.Qi, 2005 [12]):  Let $P=(p_{ij})$ be a matrix of order $n$, $\mathbb {A}$ and $\mathbb {B}$ be two order $m$ dimension $n$ tensors. If the elements of $\mathbb {A}$ and $\mathbb {B}$ satisfy the following relation:
$$b_{i_1\cdots i_m}=\sum _{j_1,\cdots ,j_m=1}^n a_{j_1\cdots j_m}p_{i_1j_1}p_{i_2j_2}\cdots p_{i_mj_m}.$$
Then it is written as $\mathbb {B}=P^m\mathbb {A}$. In the case when $P$ is a real orthogonal matrix, then $\mathbb {A}$ and $\mathbb {B}$ are called "orthogonal similar".

\vskip 0.2cm

\noindent {\bf Remark 2.3}: (1) By using the tensor product defined in section 1 and the formula \textcolor[rgb]{1.00,0.00,0.00}{(2.1)}, we see that the relation $\mathbb {B}=P^m\mathbb {A}$ defined above is actually the same as the tensor relation $\mathbb {B}=P\mathbb {A}P^T$. Since the corresponding relation for matrices is called the "congruence relation", so in this paper we choose to call the two tensors $\mathbb {A}$ and $\mathbb {B}$ to be congruent if $\mathbb {B}=P\mathbb {A}P^T$ for some matrix $P$. In the case when $P$ is a real orthogonal matrix, then we call $\mathbb {A}$ and $\mathbb {B}$ as "orthogonal congruent".

\vskip 0.18cm

\noindent (2) Notice that when $P$ is a real orthogonal matrix, the relation $P\mathbb {I}P^T=\mathbb {I}$ does not necessarily hold. That is to say, $P\mathbb {A}P^T$ and $\mathbb {A}$ are not necessary similar in the sense of our Definition 2.3. This is another reason why in this paper we choose to call them "orthogonal congruent" instead of "orthogonal similar".     \qed

\vskip 0.2cm

It is easy to see that the orthogonal congruence relation is also an equivalent relation between tensors.

\vskip 0.2cm

In [12, Theorem 7], L.Qi also proved the following result about the {\bf E}-eigenvalues and {\bf E}-eigenvectors of tensors, which can now be simply stated and proved as the following, by using our tensor product and its associative law.

\vskip 0.2cm

\noindent {\bf Theorem 2.6 }: Suppose that the two tensors $\mathbb {A}$ and $\mathbb {B}$ are  "orthogonal congruent", namely $\mathbb {B}=P\mathbb {A}P^T$ for some real orthogonal matrix $P$. Then $x$ is an
{\bf E}-eigenvector of $\mathbb {A}$ corresponding to the {\bf E}-eigenvalue $\lambda $ if and only if $y=Px$ is an {\bf E}-eigenvector of $\mathbb {B}$ corresponding to the {\bf E}-eigenvalue $\lambda $.

\vskip 0.2cm

\noindent {\bf Proof }: Since $y^Ty=x^Tx$ for $y=Px$, it is obvious that $x^Tx=1$ if and only if $y^Ty=1$. Also, by the associative law of the tensor product, we have
$$\mathbb {B}y=\lambda y\Longleftrightarrow (P\mathbb {A}P^T)(Px)=\lambda (Px)\Longleftrightarrow
P(\mathbb {A}x)=P(\lambda x)\Longleftrightarrow \mathbb {A}x=\lambda x$$

\vskip 1.88cm

\section{The direct product $\mathbb {A}\otimes \mathbb {B}$ of tensors with applications on the spectra of products of hypergraphs}

It is well-known that the direct product of matrices (denoted by $\otimes $) is a useful concept and tool in matrix theory. It has many applications in various fields, including the applications in the study of the spectra of products of graphs. In this section, we will show that the concept of the direct product of matrices can be generalized to tensors, and can also be applied to the study of the spectrum of the products of hypergraphs.

\vskip 0.18cm

\noindent {\bf Definition 3.1}: Let $\mathbb {A}$ and $\mathbb {B}$ be two order $k$ tensors with dimension $n$ and $m$, respectively. Define the direct product $\mathbb {A}\otimes \mathbb {B}$ to be the following tensor of order $k$ and dimension $nm$ (the set of subscripts is taken as $[n]\times [m]$ in the lexicographic order):
$$(\mathbb {A}\otimes \mathbb {B})_{(i_1,j_1)(i_2,j_2)\cdots (i_k,j_k)}=a_{i_1i_2\cdots i_k}b_{j_1j_2\cdots j_k} $$

\vskip 0.2cm

From this definition it is easy to see the following proposition.

\vskip 0.1cm

\noindent {\bf Proposition 3.1}:

\vskip 0.1cm

\noindent (1).  $(\mathbb {A}_1+\mathbb {A}_2)\otimes \mathbb {B}=\mathbb {A}_1\otimes \mathbb {B}+\mathbb {A}_2\otimes \mathbb {B}.$

\vskip 0.1cm

\noindent (2).  $\mathbb {A}\otimes (\mathbb {B}_1+\mathbb {B}_2)=\mathbb {A}\otimes \mathbb {B}_1+\mathbb {A}\otimes \mathbb {B}_2 .$

\vskip 0.1cm

\noindent (3).  $(\lambda \mathbb {A})\otimes \mathbb {B}=\mathbb {A}\otimes (\lambda \mathbb {B})=\lambda (\mathbb {A}\otimes \mathbb {B}).$  \qquad ($\lambda \in \mathbb {C}$.)

\vskip 0.18cm

The following theorem gives an important relation between the direct product of tensors and the general product of tensors defined in section 1, which is also a generalization of the corresponding relation for matrices.

\vskip 0.38cm

\noindent {\bf Theorem 3.1}: Let $\mathbb {A}$ and $\mathbb {B}$ be two order $k+1$ tensors with dimension $n$ and $m$, respectively. Let $\mathbb {C}$ and $\mathbb {D}$ be two order $r+1$ tensors with dimension $n$ and $m$, respectively. Then we have:
$$(\mathbb {A}\otimes \mathbb {B})(\mathbb {C}\otimes \mathbb {D})=(\mathbb {AC})\otimes (\mathbb {BD}).$$

\noindent {\bf Proof. } For the simplicity of notation, we write
$$(\alpha ,\beta )=((a_1,b_1),\cdots , (a_r,b_r))\in ([n]\times [m])^r, \quad (\mbox {if } \alpha=(a_1,\cdots ,a_r)\in [n]^r, \  \ \beta=(b_1,\cdots , b_r)\in [m]^r ). $$
In the following, we assume $\alpha_i\in [n]^r$, $\beta_j\in [m]^r$,  $ \ (i,j=1,\cdots,k)$. Then by Definition 1.1 we have
$$\begin{aligned}
\displaystyle &  \left ( (\mathbb {A}\otimes \mathbb {B})(\mathbb {C}\otimes \mathbb {D})  \right )_{(i,j)(\alpha_1 ,\beta_1 )\cdots (\alpha_k ,\beta_k )}      \\
&= \sum_{(i_1,j_1),\cdots, (i_k,j_k)\in [n]\times [m]} (\mathbb {A}\otimes \mathbb {B})_{(i,j)(i_1,j_1)\cdots  (i_k,j_k)} (\mathbb {C}\otimes \mathbb {D})_{(i_1,j_1)(\alpha_1,\beta_1)}\cdots (\mathbb {C}\otimes \mathbb {D})_{(i_k,j_k)(\alpha_k,\beta_k)}      \\
&= \sum_{i_1,\cdots ,i_k=1}^n \sum_{j_1,\cdots ,j_k=1}^m a_{ii_1\cdots i_k}b_{jj_1\cdots j_k} (c_{i_1\alpha_1}d_{j_1\beta_1})\cdots  (c_{i_k\alpha_k}d_{j_k\beta_k})        \\
&=\left ( \sum_{i_1,\cdots ,i_k=1}^n a_{ii_1\cdots i_k}c_{i_1\alpha_1}\cdots c_{i_k\alpha_k} \right )
\left ( \sum_{j_1,\cdots ,j_k=1}^m b_{jj_1\cdots j_k} d_{j_1\beta_1}\cdots d_{j_k\beta_k}  \right ) \\
&= (\mathbb {AC})_{i\alpha_1\cdots \alpha_k}(\mathbb {BD})_{j\beta_1\cdots \beta_k}   \\
&= \left ( (\mathbb {AC})\otimes (\mathbb {BD})  \right  )_{(i,j)(\alpha_1 ,\beta_1 )\cdots (\alpha_k ,\beta_k )}    \\
\end{aligned}
$$
Thus we have $(\mathbb {A}\otimes \mathbb {B})(\mathbb {C}\otimes \mathbb {D})=(\mathbb {AC})\otimes (\mathbb {BD}).$  \qed

\vskip 0.18cm

Now we consider the applications of the direct product of tensors in the study of the spectra of the products of hypergraphs.

\vskip 0.18cm

\noindent {\bf Definition 3.2} ([5], The Cartesian product of hypergraphs): Let $G$ and $H$ be two hypergraphs.
Define the Cartesian product $G\Box H$ of $G$ and $H$ as: $V(G\Box H)=V(G)\times V(H)$, and $\{(i_1,j_1),\cdots, (i_r,j_r)\}\in E(G\Box H)$ if and only if one of the following two conditions holds:

\vskip 0.1cm

\noindent (1). $i_1=\cdots =i_r$ and $\{j_1,\cdots ,j_r\}\in E(H)$.

\vskip 0.1cm

\noindent (2). $j_1=\cdots =j_r$ and $\{i_1,\cdots ,i_r\}\in E(G)$.

\vskip 0.28cm

A hypergraph $H=(V,E)$ is called $k$-uniform if every edge of $H$ is a subset of $V$ with $k$ elements. The adjacency tensor of $H$ (under certain ordering of vertices) is the order $k$ dimension $n$ tensor $\mathbb {A}_H$ with the entries ([5])
$$a_{i_1i_2\cdots i_k}=\left \{\begin{array}{cc}
                               \frac {1} {(k-1)!} &   \mbox{if }  \ \{i_1,i_2,\cdots ,i_k\} \in E(H) \\
                               0 &   \mbox{otherwise}\\
                            \end{array}
                          \right.$$
The characteristic polynomial and spectrum of a uniform hypergraph $H$ are that of its adjacency tensor.

\vskip 0.2cm

Using the direct product of tensors, the adjacency tensor of the Cartesian product of $k$-uniform hypergraphs can be expressed as in the following theorem.

\vskip 0.28cm

\noindent {\bf Theorem 3.2}: Let $\mathbb {A}$ and $\mathbb {B}$ be the adjacency tensors of a $k$-uniform hypergraph $G$ of order $n$ and a $k$-uniform hypergraph $H$ of order $m$, respectively. Then the adjacency tensor of $G\Box H$ is $\mathbb {A}\otimes \mathbb {I}_m +  \mathbb {I}_n\otimes \mathbb {B}$ (where the ordering of the vertices of $G\Box H$ is taken to be the lexicographic ordering of the elements of the set $V(G)\times V(H)$).

\noindent {\bf Proof. } Let $\mathbb {C}$ be the adjacency tensor of the hypergraph $G\Box H$. Then by definition we can check that:
$$c_{(i_1,j_1)\cdots (i_k,j_k)}=\left \{\begin{array}{cc}
                               b_{j_1\cdots j_k}&   \mbox{if }  \ i_1=i_2=\cdots =i_k  \\
                                a_{i_1\cdots i_k} & \mbox{ if } \  j_1=j_2=\cdots =j_k  \\
                               0 &   \mbox{otherwise}\\
                            \end{array}
                          \right.  \eqno {(3.1)}$$
Notice that all the diagonal entries of $\mathbb {A}$ and $\mathbb {B}$ are zero, so it follows from (3.1) that
$$c_{(i_1,j_1)\cdots (i_k,j_k)}=a_{i_1\cdots i_k}\delta _{j_1\cdots j_k}+\delta_{i_1\cdots i_k}b_{j_1\cdots j_k}$$
It is also easy to see that
$$(\mathbb {A}\otimes \mathbb {I}_m +  \mathbb {I}_n\otimes \mathbb {B})_{(i_1,j_1)\cdots (i_k,j_k)}=a_{i_1\cdots i_k}\delta _{j_1\cdots j_k}+\delta_{i_1\cdots i_k}b_{j_1\cdots j_k}$$
So we have $\mathbb {C}=\mathbb {A}\otimes \mathbb {I}_m +  \mathbb {I}_n\otimes \mathbb {B}$.
\qed

\vskip 0.18cm

Now we consider another kind of product of the hypergraphs.

\vskip 0.28cm

\noindent {\bf Definition 3.3} (The direct product of hypergraphs): Let $G$ and $H$ be two hypergraphs.
Define the direct product $G\times H$ of $G$ and $H$ as: $V(G\times H)=V(G)\times V(H)$, and $\{(i_1,j_1),\cdots, (i_r,j_r)\}\in E(G\times H)$ if and only if $\{i_1,\cdots ,i_r\}\in E(G)$ and $\{j_1,\cdots ,j_r\}\in E(H)$.

\vskip 0.28cm

Using the direct product of tensors, the adjacency tensor of the direct product of $k$-uniform hypergraphs can be easily obtained as in the following theorem.

\vskip 0.28cm

\noindent {\bf Theorem 3.3}: Let $\mathbb {A}$ and $\mathbb {B}$ be the adjacency tensors of a $k$-uniform hypergraph $G$ of order $n$ and a $k$-uniform hypergraph $H$ of order $m$, respectively. Then the adjacency tensor of $G\times H$ is $(k-1)!(\mathbb {A}\otimes \mathbb {B})$ (where the ordering of the vertices of $G\times H$ is taken to be the lexicographic ordering of the elements of the set $V(G)\times V(H)$).

\noindent {\bf Proof. } Let $\mathbb {C}$ be the adjacency tensor of the hypergraph $G\times H$. Then by definition we can check that
$$c_{(i_1,j_1)\cdots (i_k,j_k)}=\left \{\begin{array}{cc}
                               \frac {1} {(k-1)!} &   \mbox{if }  \ a_{i_1\cdots i_k} \mbox{ and } \ b_{j_1\cdots j_k} \ \mbox { are both nonzero } \\
                               0 &   \mbox{otherwise}\\
                            \end{array}
                          \right.$$
Thus it follows that
$$c_{(i_1,j_1)\cdots (i_k,j_k)}=(k-1)!a_{i_1\cdots i_k}b_{j_1\cdots j_k}=(k-1)!(\mathbb {A}\otimes \mathbb {B})_{(i_1,j_1)\cdots (i_k,j_k)}$$
So we have $\mathbb {C}=(k-1)!(\mathbb {A}\otimes \mathbb {B})$.
\qed

\vskip 0.18cm

The following theorem gives the relation between the eigenvalue-eigenvectors of the tensors $\mathbb {A}$ and $\mathbb {B}$ and that of $\mathbb {A}\otimes  \mathbb {B}$ and $\mathbb {A}\otimes \mathbb {I}_m +  \mathbb {I}_n\otimes \mathbb {B}$.

\vskip 0.38cm

\noindent {\bf Theorem 3.4}: Let $\mathbb {A}$ and $\mathbb {B}$ be two order $k$ tensors with dimension $n$ and $m$, respectively. Suppose that we have $\mathbb {A}u=\lambda u^{[k-1]}$, and $\mathbb {B}v=\mu v^{[k-1]}$, and we also write $w=u\otimes v$. Then we have:

\vskip 0.1cm

\noindent (1).  $(\mathbb {A}\otimes \mathbb {I}_m +  \mathbb {I}_n\otimes \mathbb {B})w=(\lambda +\mu)w^{[k-1]}.$

\vskip 0.1cm

\noindent (2).  $(\mathbb {A}\otimes  \mathbb {B})w=(\lambda \mu)w^{[k-1]} .$

\noindent {\bf Proof. }  We have by Proposition 3.1 and Theorem 3.1 that:

\vskip 0.1cm

\noindent (1).
$$\begin{aligned}
\displaystyle &   (\mathbb {A}\otimes \mathbb {I}_m +  \mathbb {I}_n\otimes \mathbb {B})w=(\mathbb {A}\otimes \mathbb {I}_m +  \mathbb {I}_n\otimes \mathbb {B})(u\otimes v)= (\mathbb {A}u)\otimes (\mathbb {I}_m v)+ ( \mathbb {I}_n u)\otimes (\mathbb {B}v)\\
&=(\lambda u^{[k-1]})\otimes (v^{[k-1]}) + ( u^{[k-1]})\otimes (\mu v^{[k-1]})  \\
&=(\lambda +\mu)(u^{[k-1]}\otimes v^{[k-1]})=(\lambda +\mu)(u\otimes v)^{[k-1]}= (\lambda +\mu)w^{[k-1]}\\
\end{aligned}
$$

\vskip 0.1cm

\noindent (2).
$$\begin{aligned}
\displaystyle &  (\mathbb {A}\otimes  \mathbb {B})w= (\mathbb {A}\otimes  \mathbb {B}) (u\otimes v) =   (\mathbb {A}u)\otimes  (\mathbb {B}v)= (\lambda u^{[k-1]})\otimes (\mu v^{[k-1]})\\
&=(\lambda \mu)(u^{[k-1]}\otimes v^{[k-1]})=(\lambda \mu)(u\otimes v)^{[k-1]}= (\lambda \mu)w^{[k-1]}.\\
\end{aligned}
$$
\qed

From Theorem 3.4, we can obtain the following results about the spectra of the hypergraphs $G\Box H$ and $G\times H$, where the result on the spectrum of $G\Box H$ was obtained by Cooper and Dutle in [5], but by using the direct product of tensors defined here, our Theorem 3.4 actually gives a simplified proof of that result on the spectrum of $G\Box H$.

\vskip 0.38cm

\noindent {\bf Theorem 3.5}: Let $G$ and $H$ be two $k$-uniform hypergraphs. Let $\lambda $ be an eigenvalue of $G$ with corresponding eigenvector $u$, and $\mu $ be an eigenvalue of $H$ with corresponding eigenvector $v$,       respectively. Then we have:

\vskip 0.1cm

\noindent (1). ([5]) $ \ \lambda +\mu$ is an eigenvalue of $G\Box H$ with corresponding eigenvector $u\otimes v$.

\vskip 0.1cm

\noindent (2).  $\lambda \mu$ is an eigenvalue of $G\times H$ with corresponding eigenvector $u\otimes v$.
\qed

\vskip 1.88cm

\section{Applications in the study of nonnegative tensors}

In [2] and [4], Chang et al studied the properties of the spectra of nonnegative tensors. They defined the irreducibility of tensors, and the primitivity of nonnegative tensors (also see [11]), and extended the well-known Perron-Frobenius Theorem from the nonnegative irreducible matrices to nonnegative irreducible tensors. They also extended many important properties of primitive matrices to primitive tensors. In [14] and [15], Q.Yang and Y.Yang gave some further results on the generalizations of Perron-Frobenius Theorem, including some results on the symmetry of the spectra of nonnegative irreducible tensors.

\vskip 0.18cm

In this section, we will give some applications of the tensor product in the study of nonnegative tensors, especially in the study of primitive tensors. Among these applications, we will give a simple characterization of the primitive tensors in terms of the zero patterns of the powers of $\mathbb {A}$ (in the sense of the tensor product defined in section 1), and give some upper bounds for the primitive degrees of the primitive tensors. We also give some result on the cyclic index of a nonnegative irreducible tensor.

\vskip 0.18cm

Let $Z(\mathbb {A})$ be the tensor obtained by replacing all the nonzero entries of $\mathbb {A}$ by one. Then $Z(\mathbb {A})$ is called the zero-nonzero pattern of $\mathbb {A}$ (or simply the zero pattern of $\mathbb {A}$).

\vskip 0.1cm

First we have the following result on the zero patterns of the products of two nonnegative tensors.

\vskip 0.18cm

\noindent {\bf Lemma 4.1}: The zero pattern of the product $\mathbb {AB}$ of two nonnegative tensors $\mathbb {A}$ and $\mathbb {B}$ is uniquely determined by the zero patterns of $\mathbb {A}$ and $\mathbb {B}$.

\noindent {\bf Proof. }  Suppose that $Z(\mathbb {A}_1)=Z(\mathbb {A}_2)$ and $Z(\mathbb {B}_1)=Z(\mathbb {B}_2)$. Then there exist positive numbers $c_1,c_2,M_1$ and $M_2$ such that
$$c_1\mathbb {A}_1\le \mathbb {A}_2\le M_1\mathbb {A}_1; \qquad \quad  c_2\mathbb {B}_1\le \mathbb {B}_2\le M_2\mathbb {B}_1$$
From this it follows that $c_1c_2^{m-1}\mathbb {A}_1\mathbb {B}_1\le \mathbb {A}_2\mathbb {B}_2\le M_1M_2^{m-1}\mathbb {A}_1\mathbb {B}_1$ (where $m$ is the order of $\mathbb {A}_1$). Thus we have
$Z(\mathbb {A}_1\mathbb {B}_1)=Z(\mathbb {A}_2\mathbb {B}_2)$.  \qed

\vskip 0.1cm

The following result will be used in Definition 4.1.

\vskip 0.18cm

\noindent {\bf Proposition 4.1}: Let $\mathbb {A}$ be an order $m$ dimension $n$ nonnegative tensor. Then the following three conditions are equivalent:

\vskip 0.1cm

\noindent (1). For any $i,j\in [n]$, $ \ a_{ij\cdots j}>0$ holds.

\vskip 0.1cm

\noindent (2). For any $j\in [n]$, $ \ \mathbb {A}e_j>0$ holds (where $e_j$ is the $j^{th}$ column of the identity matrix $I_n$).

\vskip 0.1cm

\noindent (3). For any nonnegative nonzero vector $x\in \mathbb {R}^n$, $ \ \mathbb {A}x>0$ holds.

\vskip 0.2cm

\noindent The proof of this proposition is straightforward. \qed

\vskip 0.18cm

\noindent {\bf Definition 4.1} (Pearson, [10]): A nonnegative tensor $\mathbb {A}$ is called essentially positive, if it satisfies one of the three conditions in Proposition 4.1.

\vskip 0.18cm

In [4] and [11], Chang et al and Pearson define the primitive tensors as follows.

\vskip 0.2cm

\noindent {\bf Definition 4.2} ([4] and [11]): Let $\mathbb {A}$ be a nonnegative tensor with order $m$ and dimension $n$. Define the map $T_{\mathbb {A}}$ from $\mathbb {R}^n$ to $\mathbb {R}^n$ as: $T_{\mathbb {A}}(x)=(\mathbb {A}x)^{[\frac {1} {m-1}]}$ (here $\mathbb {A}x$ is understood as the tensor product defined in section 1). If there exists some positive integer $r$ such that $T_{\mathbb {A}}^r(x)>0$ for all nonnegative nonzero vectors $x\in \mathbb {R}^n$, then $\mathbb {A}$ is called primitive, and the smallest such integer $r$ is called the primitive degree of $\mathbb {A}$, denoted by $\gamma (\mathbb {A})$.

\vskip 0.1cm

By using the properties of tensor product and Lemma 4.1, we can obtain the following simple characterization for primitive tensors.

\vskip 0.38cm

\noindent {\bf Theorem 4.1}: A nonnegative tensor $\mathbb {A}$ is primitive if and only if there exists some positive integer $r$ such that $\mathbb {A}^r$ is essentially positive. Furthermore, the smallest such $r$ is the primitive degree of $\mathbb {A}$.

\vskip 0.18cm

\noindent {\bf Proof}: For convenience, here we use the notation $\mathbb {A}\backsimeq \mathbb {B}$ to denote that $\mathbb {A}$ and $\mathbb {B}$ have the same zero pattern.

\vskip 0.18cm

We first use induction on $r$ to show that for any nonnegative nonzero vectors $x\in \mathbb {R}^n$, $T_{\mathbb {A}}^r(x)$ and $\mathbb {A}^rx$ have the same zero pattern. When $r=1$, we obviously have $T_{\mathbb {A}}(x)=(\mathbb {A}x)^{[\frac {1} {m-1}]}\backsimeq \mathbb {A}x$. In general, by the inductive hypothesis, Lemma 4.1 and the associative law of the tensor product we have:
$$T_{\mathbb {A}}^r(x)=T_{\mathbb {A}}(T_{\mathbb {A}}^{r-1}(x))=(\mathbb {A}T_{\mathbb {A}}^{r-1}(x))^{[\frac {1} {m-1}]}\backsimeq \mathbb {A}T_{\mathbb {A}}^{r-1}(x)\backsimeq \mathbb {A}(\mathbb {A}^{r-1}x)
=\mathbb {A}^rx.$$

From this we have for any nonzero vector $x\ge 0$ and positive integer $r$ that:
$T_{\mathbb {A}}^r(x)>0\Longleftrightarrow \mathbb {A}^rx>0$. Thus our result follows directly.
\qed

\vskip 0.18cm

Using the similar ideas, the characterization of (nonnegative) irreducible tensors given in [15, Theorem 5.2] can also be simply stated as: A nonnegative tensor $\mathbb {A}$ of dimension $n$ is irreducible if and only if $(\mathbb {A}+\mathbb {I})^{n-1}$ is essentially positive.

\vskip 0.28cm

The following Lemma 4.2 will be used in the proof of Theorem 4.2.

\vskip 0.38cm

\noindent {\bf Lemma 4.2}: Let $\mathbb {A}$ be a nonnegative tensor. If $\mathbb {A}^r$ is essentially positive, then $\mathbb {A}^{r+1}$ is also essentially positive.

\vskip 0.18cm

\noindent {\bf Proof}: By hypothesis, $\mathbb {A}^rx>0$ for all nonzero vectors $x\ge 0$. It follows that $y=\mathbb {A}x$ is also a nonzero (nonnegative) vector (otherwise we will have $\mathbb {A}^rx=0$). Therefore we have (for any nonzero vectors $x\ge 0$) that: $\mathbb {A}^{r+1}x=\mathbb {A}^r(\mathbb {A}x)=\mathbb {A}^ry>0$, so $\mathbb {A}^{r+1}$ is  essentially positive.
\qed

\vskip 0.18cm

The following theorem shows that, although there are infinitely many primitive tensors of order $m$ and dimension $n$, their primitive degrees are bounded above and have the following upper bound.

\vskip 0.38cm

\noindent {\bf Theorem 4.2} (An upper bound of the primitive degrees): Let $\mathbb {A}$ be a nonnegative primitive tensor with order $m$ and dimension $n$. Then its primitive degree $\gamma (\mathbb {A})\le 2^{n^m}$.

\vskip 0.18cm

\noindent {\bf Proof}: First we notice that there are altogether $2^{n^m}$ many different zero patterns for order $m$ dimension $n$ tensors. So among the powers $\mathbb {A}, \cdots , \mathbb {A}^{2^{n^m}+1}$, there must be some two powers having the same zero pattern. Thus there exist some integer $k\le 2^{n^m}$ and some integer $p>0$ such that $Z(\mathbb {A}^k)= Z(\mathbb {A}^{k+p})$. From this it will follow by Lemma 4.1 that $Z(\mathbb {A}^k)= Z( \mathbb {A}^{k+tp})$ for all nonnegative integers $t$. On the other hand, by Lemma 4.2 we also know that $\mathbb {A}^{k+tp}$ is essentially positive for all sufficiently large $t$, so $\mathbb {A}^k$ is essentially positive, which implies by definition that $\gamma (\mathbb {A})\le k\le 2^{n^m}$.
\qed

\vskip 0.28cm

\noindent {\bf Remark 4.1}: The upper bound $2^{n^m}$ can be replaced by the number of different zero patterns of all the primitive tensors of order $m$ and dimension $n$. Thus the upper bound in Theorem 4.2 is not sharp.  \qed

\vskip 0.58cm

In [10], Pearson defined the majorization matrix $M(\mathbb {A})$ of a tensor $\mathbb {A}$ as the following:
$$(M(\mathbb {A}))_{ij}=a_{ij\cdots j} \qquad \quad  (i,j=1,\cdots ,n)$$
Pearson ([10]) also proved that if $M(\mathbb {A})$ is an irreducible matrix, then $\mathbb {A}$ is also irreducible.

Also, the spectral radius of $\mathbb {A}$ is denoted by $\rho (\mathbb {A})$, and the cyclic index of a nonnegative irreducible tensor $\mathbb {A}$ is defined to be the number of distinct eigenvalues of $\mathbb {A}$ whose absolute value is $\rho (\mathbb {A})$.

\vskip 0.2cm

In the following theorem, we use the majorization matrix $M(\mathbb {A})$ to study the
cyclic index of a nonnegative tensor $\mathbb {A}$.

\vskip 0.28cm

\noindent {\bf Theorem 4.3}: Let $\mathbb {A}$ be a nonnegative tensor with order $m$ and dimension $n$. If the matrix $M(\mathbb {A})$ is irreducible, then the cyclic index of $M(\mathbb {A})$ is a multiple of the cyclic index of $\mathbb {A}$.

\vskip 0.18cm

\noindent {\bf Proof}: By ([10]) we know that $\mathbb {A}$ is also irreducible. Now let $k$ be the cyclic index of $\mathbb {A}$, and $p$ be the cyclic index of the irreducible matrix $M(\mathbb {A})$. If $k=1$, the result holds obviously. So we assume that $k\ge 2$.

Let $\lambda =\rho (\mathbb {A})e^{i2\pi /k}$, let $y=(y_1,\cdots ,y_n)^T$ be an eigenvector of $\mathbb {A}$ (with full support) corresponding to the eigenvalue $\lambda $, and write (as in [4])
$$y_j=|y_j|e^{i2\pi \phi_j}\qquad \quad  (j=1,\cdots ,n)$$
Then by [4, Lemma 4.4], we have for any $a_{i_1\cdots i_m}>0$ that
$$\phi_{i_2}+\cdots +\phi_{i_m}-(m-1)\phi_{i_1}-\frac {1} {k}\in \mathbb {Z}\qquad (\mathbb {Z} \  \  \mbox {is the set of all integers} )\eqno {(4.1)}$$

Let $j_1\cdots j_rj_1$ be a (directed) cycle of the associated digraph $D(M(\mathbb {A}))$ ([1]) of the matrix $M(\mathbb {A})$ (with length $r$), then we have
$$\begin{aligned}
\displaystyle & (M(\mathbb {A}))_{j_tj_{t+1}}>0 \  \ (t=1,\cdots ,r, \ mod \ r)\Longrightarrow a_{j_tj_{t+1}\cdots j_{t+1}}>0 \  \ (t=1,\cdots ,r, \ mod \ r)  \\
& \Longrightarrow (m-1)\phi_{j_{t+1}}-(m-1)\phi_{j_t}-\frac {1} {k}\in \mathbb {Z} \quad  \ (t=1,\cdots ,r, \ mod \ r) \\
\end{aligned}
$$
Adding the above $r$ relations together, we obtain
$$(m-1)\sum_{t=1}^r(\phi_{j_{t+1}}-\phi_{j_t})-\frac {r} {k}=-\frac {r} {k} \in \mathbb {Z} \eqno {(4.2)}$$
Thus $k$ divides $r$. But $r$ is an arbitrary cycle length of the associated digraph of the matrix $M(\mathbb {A})$, so $k$ divides the greatest common divisor of all the cycle lengths of the associated digraph of $M(\mathbb {A})$, which is just $p$, the cyclic index of the irreducible matrix $M(\mathbb {A})$ ([1]).
\qed

\vskip 0.1cm

We also have the following relation between the primitive degrees of $\mathbb {A}$ and $M(\mathbb {A})$.

\vskip 0.38cm

\noindent {\bf Corollary 4.1}: Let $\mathbb {A}$ be a nonnegative tensor with order $m$ and dimension $n$. If $M(\mathbb {A})$ is primitive, then $\mathbb {A}$ is also primitive and in this case, we have
$\gamma (\mathbb {A})\le \gamma (M(\mathbb {A}))\le (n-1)^2+1$.

\vskip 0.18cm

\noindent {\bf Proof}: $\gamma (M(\mathbb {A}))\le (n-1)^2+1$ is the well-known Wielandt's upper bound for primitive matrices ([16]). So we only need to show that $\mathbb {A}$ is primitive and $\gamma (\mathbb {A})\le \gamma (M(\mathbb {A}))$. For this purpose, we use the notation
$\mathbb {A}\succsim \mathbb {B}$ to denote (the zero patterns) $Z(\mathbb {A})\ge Z(\mathbb {B})$. Then we  first prove for two nonnegative tensors $\mathbb {A}$ and $\mathbb {B}$ that $M(\mathbb {AB})\succsim M(\mathbb {A})M(\mathbb {B})$. For if $(M(\mathbb {A})M(\mathbb {B}))_{ij}>0$, then there exists $h\in [n]$ such that $(M(\mathbb {A}))_{ih}>0$ and $(M(\mathbb {B}))_{hj}>0$. It means that $a_{ih\cdots h}>0$ and $b_{hj\cdots j}>0$, which implies that
$$(M(\mathbb {AB}))_{ij}=(\mathbb {AB})_{ij\cdots j}=\sum_{i_2,\cdots ,i_m=1}^n a_{ii_2\cdots i_m}b_{i_2j\cdots j}\cdots b_{i_mj\cdots j}\ge a_{ih\cdots h}(b_{hj\cdots j})^{m-1}>0.$$

Now from $M(\mathbb {AB})\succsim M(\mathbb {A})M(\mathbb {B})$ we can see that for any integer $k>0$, we have $M(\mathbb {A}^k)\succsim (M(\mathbb {A}))^k$. Take $k=\gamma (M(\mathbb {A}))$, then we have $M(\mathbb {A}^k)\succsim (M(\mathbb {A}))^k>0$, which means that $\mathbb {A}^k$ is essentially positive. So $\mathbb {A}$ is also primitive and $\gamma (\mathbb {A})\le k=\gamma (M(\mathbb {A}))$. \qed

\vskip 0.28cm

Now we propose the following conjecture on the primitive degrees.

\vskip 0.28cm

\noindent {\bf Conjecture 1}: When $m$ is fixed, then there exists some polynomial $f(n)$ on $n$ such that $\gamma (\mathbb {A})\le f(n)$ for all nonnegative primitive tensors of order $m$ and dimension $n$.

\vskip 0.2cm

In the case of $m=2$ ($\mathbb {A}$ is a matrix), the Wielandt's upper bound tells us that we can take $f(n)=(n-1)^2+1$.

\vskip 0.3cm

\noindent {\bf Definition 4.3}: Let $\mathbb {A}$ be a nonnegative tensor with order $m$ and dimension $n$. If
there exists some positive integer $k$ such that $\mathbb {A}^k>0$ is a positive tensor, then $\mathbb {A}$ is called strongly primitive, and the smallest such $k$ is called the strongly primitive degree of $\mathbb {A}$.

\vskip 0.18cm

It is obvious that the strong primitivity implies the primitivity. Also, in the matrix case ($m=2$), it is well-known ([1]) that a nonnegative matrix $A$ is primitive if and only if it is strongly primitive. But the following example shows that these two concepts are not equivalent in the case $m\ge 3$.

\vskip 0.2cm

\noindent {\bf Example 4.1}: Let $\mathbb {A}$ be a nonnegative tensor of order $m\ge 3$ and dimension $n$, and   with
$$a_{i_1i_2\cdots i_m}=\left\{
                            \begin{array}{cc}
                              1 &   \mbox{if } i_2=\cdots =i_m  \\
                              0 &  \mbox{otherwise}\\
                            \end{array}
                          \right.$$
Then $\mathbb {A}$ is a primitive tensor, but not a strongly primitive tensor.

\vskip 0.2cm

\noindent {\bf Proof}: $\mathbb {A}$ is primitive since $\mathbb {A}$ itself is already essentially positive. Now we use induction on $k$ to show that if $i_2\ne i_3$, then $(\mathbb {A}^k)_{i_1i_2i_3\cdots i_r}=0$ (where $r=(m-1)^k+1$).

If $k=1$, the result follows from the definition of $\mathbb {A}$. In general when $k\ge 2$, let $\mathbb {B}=\mathbb {A}^{k-1}=(b_{j_1\cdots j_t})$, where $t=(m-1)^{k-1}+1$. Then by the inductive hypothesis, we have $b_{j_1i_2i_3\cdots j_t}=0$ when $i_2\ne i_3$. Therefore we have
$$(\mathbb {A}^k)_{i_1i_2i_3\cdots i_r}=(\mathbb {AB})_{i_1i_2i_3\cdots i_r}=\sum_{j_1,\cdots ,j_{m-1}=1}^n a_{i_1j_1\cdots j_{m-1}}b_{j_1i_2i_3\cdots i_t}b_{j_2i_{t+1}\cdots i_{2t-1}}\cdots b_{j_{m-1}i_{r-t+2}\cdots i_r}=0 \quad (i_2\ne i_3)$$
Thus $\mathbb {A}^k$ is not a positive tensor for all integers $k$, so $\mathbb {A}$ is not strongly primitive.
\qed

\vskip 1.8cm
\noindent{\bf References}

\bibliographystyle{acm}

\begin{thebibliography}{10}

\bibitem{brualdi}
{\sc R.A.Brualdi and H.J.Ryser,}
\newblock Combinatorial matrix theory.
\newblock {\em Cambridge University Press, 1991.\/}



\bibitem{chang2008}
{\sc K.C.Chang, K.Pearson, and T.Zhang,}
\newblock Perron-Frobenius theorem for nonnegative tensors.
\newblock {\em Commun. Math. Sci. 6(2008), 507-520.\/}


\bibitem{chang}
{\sc K.C.Chang, K.Pearson, and T.Zhang,}
\newblock On eigenvalues of real symmetric tensors.
\newblock {\em J. Math. Anal. Appl. 350(2009), 416-422.\/}


\bibitem{chang2011}
{\sc K.C.Chang, K.Pearson, and T.Zhang,}
\newblock Primitivity, the convergence of the NQZ method, and the largest eigenvalue for nonnegative tensors.
\newblock {\em SIAM J. Matrix Anal. Appl., 32(2011), 806-819.\/}


\bibitem{Cooper}
{\sc J.Cooper and A.Dutle,}
\newblock {Spectra of uniform hypergraphs}.
\newblock {\em Linear Alg. Appl., 436(2012), 3268-3292.\/}

\bibitem{Cox}
{\sc D.Cox, J.Little and D.O'Shea,}
\newblock {Using algebraic geometry}.
\newblock {\em Springer-Verlag, New York, 1998.\/}

\bibitem{Gelfand}
{\sc I.Gelfand, M.Kapranov and A.Zelevinsky,}
\newblock Discriminants, Resultants and Multidimensional Determinants
\newblock {\em Birkhauser Boston, 1994.\/}

\bibitem{lim}
{\sc L.H.Lim,}
\newblock Singular values and eigenvalues of tensors, a variational approach,
\newblock {\em in Proceedings 1st IEEE international workshop on computational advances of multitensor adaptive processing (2005), 129-132.\/}


\bibitem{Ng}
{\sc M.Ng, L.Qi and G.Zhou,}
\newblock Finding the largest eigenvalue of a nonnegative tensor.
\newblock {\em SIAM J. Matrix Anal. Appl., 31(2009), 1090-1099.\/}


\bibitem{Pearson2010}
{\sc K.Pearson,}
\newblock {\em Essentially positive tensors}.
\newblock {\em Int.J.Algebra, 4(2010), 421-427.\/}

\bibitem{Pearson}
{\sc K.Pearson,}
\newblock {Primitive tensors and convergence of an iterative process for the eigenvalues of a primitive tensor}.
\newblock {\em arXiv: 1004-2423v1, 2010.\/}

\bibitem{Qi2005}
{\sc L.Qi,}
\newblock Eigenvalues of a real supersymmetric tensor.
\newblock {\em Symbolic Comput., 40(2005), 1302-1324.\/}

\bibitem{Qi-Sun}
{\sc L.Qi, W.Sun and Y.Wang,}
\newblock Numerical multilinear algebra and its applications.
\newblock {\em Front. Math. China, 2(2007), 501-526\/}


\bibitem{yy2010}
{\sc Y.Yang and Q.Yang,}
\newblock Further results for Perron-Frobenius theorem for nonnegative tensors.
\newblock {\em SIAM J. Matrix Anal. Appl., 31(2010), 2517-2530.\/}

\bibitem{yy2011}
{\sc Y.Yang and Q.Yang,}
\newblock {Further results for Perron-Frobenius theorem for nonnegative tensors II}.
\newblock {\em SIAM J. Matrix Anal. Appl., 32(2011), 1236-1250.\/}

\bibitem{wielandt}
{\sc H.Wielandt,}
\newblock Unzerlegbare, nicht negative matrizen.
\newblock {\em Math. Zeit, 52(1950), 642-648.\/}



\end{thebibliography}

\end{document}